\makeatletter

\def\eqalign#1{\,\vcenter{\openup\jot\m@th
  \ialign{\strut\hfil$\displaystyle{##}$&$\displaystyle{{}##}$\hfil
      \crcr#1\crcr}}\,}
\def\eqalignno#1{\displ@y \tabskip\@centering
  \halign to\displaywidth{\hfil$\displaystyle{##}$\tabskip\z@skip
    &$\displaystyle{{}##}$\hfil\tabskip\@centering
    &\llap{$##$}\tabskip\z@skip\crcr
    #1\crcr}}
\def\leqalignno#1{\displ@y \tabskip\@centering
  \halign to\displaywidth{\hfil$\displaystyle{##}$\tabskip\z@skip
    &$\displaystyle{{}##}$\hfil\tabskip\@centering
    &\kern-\displaywidth\rlap{$##$}\tabskip\displaywidth\crcr
    #1\crcr}}

\makeatother

\newdimen\pixel \pixel=.00333333 true in

\documentclass[11pt]{article}
\usepackage{fullpage,epsfig}  

\makeatletter
\def\bigpar{\bigbreak\@afterindentfalse\@afterheading\ignorespaces}
\def\medpar{\medbreak\@afterindentfalse\@afterheading\ignorespaces}
\def\smallpar{\smallbreak\@afterindentfalse\@afterheading\ignorespaces}

\newlength{\saveindent}
\newenvironment{proof}%
      {\bigpar{\bf Proof:}\ 
             \setlength{\saveindent}{\parindent} 
                       \ignorespaces}%
      {\stopproof\ignorespaces\bigbreak \setlength{\parindent}{\saveindent}}

      {\bigpar{\bf Proof:}%
             \setlength{\saveindent}{\parindent} 
                       \,\ignorespaces}%
      {\ignorespaces\bigbreak \setlength{\parindent}{\saveindent}}

      {\bigpar{\bf Proof:}\ %
             \setlength{\saveindent}{\parindent} 
                       \ignorespaces}%
      {\ignorespaces\bigbreak \setlength{\parindent}{\saveindent}}

\newenvironment{proofof}[1]%
      {\bigpar{\bf#1:}\ %
             \setlength{\saveindent}{\parindent} 
                       \ignorespaces}%
      {\stopproof\ignorespaces\bigbreak \setlength{\parindent}{\saveindent}}

      {\smallpar{\bf Remark:}\ 
                       \ignorespaces}%
      {\stopproof\ignorespaces\medbreak \setlength{\parindent}{\saveindent}} 

\newenvironment{remark*}%
      {\smallpar{\bf Remark:}\ 
                       \ignorespaces}%
      {\ignorespaces\medbreak \setlength{\parindent}{\saveindent}} 

      {\smallpar{\bf Remarks:}\ 
                       \ignorespaces}%
      {\stopproof\ignorespaces\medbreak \setlength{\parindent}{\saveindent}}

\newenvironment{remarks*}%
      {\smallpar{\bf Remarks:}\ 
                       \ignorespaces}%
      {\ignorespaces\medbreak \setlength{\parindent}{\saveindent}}

      {\smallpar{\bf Remark:}\ %
                       \ignorespaces}%
      {\stopproof\ignorespaces\medbreak \setlength{\parindent}{\saveindent}}

      {\smallpar{\bf Remarks:}\ %
                       \ignorespaces}%
      {\stopproof\ignorespaces\medbreak \setlength{\parindent}{\saveindent}}
\makeatother

\newtheorem{theorem}{Theorem}
\newtheorem{lemma}[theorem]{Lemma}

\newtheorem{proposition}[theorem]{Proposition}

\newtheorem{example}{Example}

\def\begex{\begin{example}\parindent=0pt \rm}
\def\endex{\end{example}}
\def\square{\vbox{\hrule height.2pt\hbox{\vrule width.2pt height5pt \kern5pt
                                   \vrule width.2pt} \hrule height.2pt}}
\def\stopproof{\hfill \square \smallskip}

\def \lphi {{u}}

\def \k {{K}}
\def \kh {{ \widehat{K}}}

\def\half{{\textstyle{1\over2}}}

\def\eighth{{\textstyle{1\over8}}}

\def\lc {{\overline L}}

\def\Sc {{\overline{S}}}

\def\zc {{\overline Z}}
\def \kh {{ \widehat{K}}}

\def \lp {{ \overline{p}}}
\def \lP {{ \overline{P}}}
\def \lphi {{ \overline{\Phi}}}
\def \lpsi {{ \overline{\psi}}}

\def\hc {{\overline{h}}}

\def\muxn{{\mu_n^x}}

\def\xt{{\overline{X}}}

\def \d {\partial}
\def\r|{{\Bigr\vert}}
\def\l|{{\Bigl\vert}}

\def \R {{\bf R}}

\def \lov {{Lov\'asz }}
\def\phi {\Phi}
\def\psistar {\psi_*}

\def\e{\epsilon}

\def\tmix{\tau}

\def\varepsilon{\mathchar"122 }

\def \one {{\mathbf 1}}

\def\p {{ \bf P}}
\def\P {{ \bf P}}

\def\eh {{ \bf {\widehat E}}}
\def\e {{ \bf {E}}}

\def\Sq{{\cal S}_q}
\def\Sq-{{\cal S}_{q-1}}

\def \lam {{\Lambda}}

\def\Px{\P_{\! \{x\}} \!}
\def\Ex{\e_{\{x\}}}
\def\wS{\widetilde{S}}
\def\hp  {   \overleftarrow{p}           \,}
\def\pim{\pi_*}

\def\one{{\mathbf 1}}
\def\RR{R}
\def\f{f}

\newcommand{\be}{\begin{equation}}
\newcommand{\ee}{\end{equation}}
\newcommand{\lab}{\label}
\begin{document}

\title{Evolving sets, mixing and heat kernel bounds}
\author{
{\sc Ben Morris}\thanks{Department of Statistics,
Evans Hall, University of California, Berkeley CA~94720.  
\newline Email:
{\tt morris@stat.berkeley.edu}.  Supported by NSF 
post-doctoral fellowship.}
\and
{\sc Yuval Peres}\thanks{Departments of Statistics and Mathematics,
 University of California, Berkeley. \newline Email:
{\tt peres@stat.berkeley.edu} Research supported 
in part by NSF Grants \#DMS-0104073 and\#DMS-0244479. . }
}
\maketitle
\thispagestyle{empty}
\begin{abstract}
We show that a new probabilistic technique,
recently introduced by the first author,
yields the sharpest bounds obtained to date on mixing times 
of Markov chains in
terms of isoperimetric properties of the state space
(also known as conductance bounds or Cheeger inequalities). 
We prove that the bounds for mixing time in total variation
obtained by Lov\'asz and Kannan, can be refined to apply to the 
maximum relative deviation $|p^n(x,y)/\pi(y) -1 |$
of the distribution at time $n$ from the stationary distribution
$\pi$.   
We then extend our results to Markov chains on infinite state spaces 
and to continuous-time chains.
Our approach yields a direct link between isoperimetric 
inequalities and heat kernel bounds;
previously, this link rested on analytic estimates
known as Nash inequalities.
\end{abstract}

\section{Introduction} \lab{intro}
\noindent
It is well known that the absence of ``bottlenecks'' in the state
space of a Markov chain implies rapid mixing.
Precise formulations of this principle, related to
Cheeger's inequality in differential geometry, have been proved by
algebraic and combinatorial techniques
\cite{Al,LS,JS,Mi,Fi,LK}. They have been used to 
approximate permanents, to sample
from the lattice points in a convex set, to estimate volumes,
and to analyze a random walk on a percolation cluster in a box.

In this paper, we show that a new probabilistic technique,
introduced in \cite{Mo},
yields the sharpest bounds obtained to date on mixing times in
terms of bottlenecks.

Let $\{p(x,y)\}$ be transition probabilities
for an irreducible Markov chain on a countable state space $V$, 
with stationary distribution $\pi$
(i.e., $\sum_{x \in V} \pi(x) p(x,y)= \pi(y)$ for all $x \in V$).
For $x, y \in V$, let $Q(x,y) = \pi(x) p(x,y)$, and 
for $S,A\subset V$, define $Q(S,A) = \sum_{s \in S, a \in A}
 Q(s,a)$. For $S \subset V$, the ``boundary size'' of $S$ is 
measured by $|\d S| = Q(S, S^c)$. Following \cite{JS}, we call 
$\phi_S:=\frac{|\d S|}{\pi(S)}$ the {\em conductance\/} of $S$.
Write $\pim:=\min_{x \in V} \pi(x)$ and
define $\phi(r)$ for $r \in [\pim,1/2]$ by
\be \lab{defphi}
\phi(r) = \inf 
\left\{  \phi_S : \pi(S) \leq r \right\}
\, .
\ee
For $r>1/2$, let $\phi(r)=\phi_*=\phi(1/2)$. 
Define the {\it $\epsilon$-uniform mixing time}
by $$
 \tmix(\epsilon) = 
 \min\Bigl\{n: \l| \frac{p^n(x,y) - \pi(y)}{\pi(y)} \r| \leq \epsilon
 \;\;\forall \,x,y \in V               \Bigr\}.$$

Jerrum and Sinclair \cite{JS} considered chains 
that are {\em reversible\/} (i.e., $Q(x,y)=Q(y,x)$ for all $x,y \in V$)
and also satisfy
\be \lab{pain}
p(x,x) \geq 1/2 \mbox{ \rm for all } x \in V  \,.
\ee
They estimated the second
eigenvalue of $p(\cdot, \cdot)$ in terms of conductance, and derived
the bound
\be \lab{JS}
\tau(\epsilon) \leq 2\phi^{-2}_* \left(\log \frac{1}{\pi_*} + \log
  \frac{1}{\epsilon} \right) \,.
\ee
Algorithmic applications of (\ref{JS}) are described in \cite{Si}. 
Extensions of (\ref{JS}) to non-reversible chains were obtained 
by Mihail \cite{Mi} and Fill~\cite{Fi}. 
A striking new idea was introduced by \lov and Kannan~\cite{LK}, who 
realized that in geometric examples, small sets often have larger
conductance, and discovered a way to exploit this. 
Let $\Vert\mu-\nu\Vert=\half\sum_{y\in V}|\mu(y)-\nu(y)|$
be the total variation distance, and denote by
\be \lab{tauv}
\tau_V(\epsilon):=
\min\Bigl\{n:\Vert p^{n}(x, \cdot)-\pi\Vert\le \epsilon \mbox{ for all
  $x \in V$}
\Bigr\}
\ee
the $\epsilon$-mixing time in total variation.
(This can be considerably smaller than the uniform mixing time 
$\tau(\epsilon)$,
see the lamplighter walk discussed at the end of this section,
or \S \ref{con},~Remark 1.)
For reversible chains that satisfy (\ref{pain}),  
\lov and Kannan proved that
\be \lab{LK}
 \tau_V(1/4) \leq 2000 \int_{\pim}^{3/4} \frac{du}{u \phi^2(u)},
\ee
This formula
was the impetus for the present paper. Related formulae 
for infinite Markov chains were
obtained earlier from Nash inequalities and are discussed below.
(As noted in \cite{MS}, there was a small error in \cite{LK}; 
the statement above is obtained from \S 3 in the survey by 
Kannan \cite{Ka}.) 

 Note that in general, $\tau_V(\epsilon) \le \tau_V(1/4) \log_2(1/\epsilon)$.
Therefore,  ignoring constant factors,
the bound in (\ref{LK}) is tighter than the bound of
(\ref{JS}), but at the cost of  employing a weaker notion of mixing.
 
\smallskip

Our main result sharpens (\ref{LK}) to a bound on
the uniform mixing time. See Theorem \ref{hk2} for a version that 
relaxes the assumption (\ref{pain}).
We use the notation $\alpha \wedge \beta := \min\{\alpha,\beta\}$.
\begin{theorem} 
\label{hk}
 Assume (\ref{pain}). Then the $\epsilon$-uniform mixing time satisfies
\be \lab{taubound}
\tau(\epsilon) \leq 1 + \int_{ 4\pi_* }^{4/\epsilon}
\frac{4 du}{u \phi^2(u)} \,.
\ee
More precisely, 
if
\be \lab{inteq}
n \ge 1 +
\int_{ 4(\pi(x) \wedge \pi(y)) }^{4/\epsilon}
\frac{4 du}{u \phi^2(u)} \, ,
\ee
then 
\be \lab{unif}
\l| \frac{p^n(x,y) - \pi(y)}{\pi(y)} \r| \leq \epsilon.
\ee
\end{theorem}
(Recall that $\phi(r)$ is constant for $r \ge \half$.)
This result has several advantages over (\ref{LK}):
\begin{itemize}
\item The uniformity in (\ref{taubound}).
\item It yields a better bound when the approximation parameter
  $\epsilon$ is small.
\item It applies to non-reversible chains.
\item It yields an improvement of the upper bound
      on the time to achieve (\ref{unif}) when 
    $\pi(x), \pi(y)$ are larger than $\pim$. 
\item The improved constant factors make the bound (\ref{taubound})
    potentially applicable as a stopping time in simulations.
    Under a convexity condition, these factors can be improved
    further; see \S \ref{con}, Remark 3.
\end{itemize}

\smallskip

Other ways to measure bottlenecks can yield sharper
bounds. One approach, based on ``blocking conductance functions''
and restricted to the mixing time in total variation
$\tau_V$, is presented in \cite[Theorem 3]{Ka}. 
\smallskip

Another boundary gauge $\psi$ 
is defined in \S \ref{fur} of the
present paper.
For the $n$-dimensional unit hypercube, this gauge
(applied to the right class of sets, see \S \ref{con}) 
gives a bound of the right order 
$\tau(1/e) =O(n \log n)$ for the uniform mixing time. 
Previous methods of measuring bottlenecks
did not yield the right order of magnitude
for the uniform mixing time in this benchmark example.
\smallskip

Theorem \ref{hk} is related to another line of research, namely
the derivation of heat kernel estimates for Markov chains using
Nash and Sobolev inequalities. For finite Markov chains,
such estimates were obtained by Chung and Yau~\cite{CY},
and by Diaconis and Saloff-Coste~\cite{DS}.
In particular, for the special case where
$\phi$ is a power law,  the conclusion of Theorem \ref{hk}
can be obtained by combining Theorems 2.3.1 and 3.3.11 of 
Saloff Coste~\cite{Sa}. For infinite Markov chains, Nash inequalities
have been developed for general isoperimetric profiles;
see Varopoulos \cite{V},
the survey by Pittet and Saloff Coste~\cite{PS},
the book~\cite{Wo}, and especially the work of Coulhon~\cite{Co,CGP}.
Even in this highly developed subject, our
probabilistic technique yields improved estimates
when the stationary measure is not uniform. 
Suppose that $\pi$ is an infinite stationary measure on $V$ for the
transition kernel $p$.
As before, we define 
\[
Q(x,y) = \pi(x) p(x,y);  \;\;\;\;\;\;\;\;\;\; 
| \d S | = Q(S, S^c); \;\;\;\;\;\;\;\;\;\;
\phi_S:=\frac{|\d S|}{\pi(S)}.
\] 
Define $\phi(r)$ for $r \in [\pim,\infty)$ by
\be \lab{i:defphi}
\phi(r) = \inf 
\left\{  \phi_S : \pi(S) \leq r \right\}
\, .
\ee

\begin{theorem} {\bf(infinite stationary measure case)}
\label{hki}

Suppose that $0<\gamma \le \half$ and $p(x,x) \ge \gamma$ for all $x \in V$.
If
\be \lab{i:inteq2}
n \ge
1+\frac{(1-\gamma)^2}{\gamma^2} \int_{4( \pi(x) \wedge \pi(y)) }^{4/\epsilon}
\frac{4 du}{u \phi^2(u)} \,,
\ee
then  
\be 
\l| \frac{p^n(x,y)}{\pi(y)} \r| \leq \epsilon.
\ee

\end{theorem}
This Theorem is proved in Section 6. For the rest of the introduction,
we focus on the case of finite stationary measure.

\begin{figure}
\label{onestep}
\centering
\epsfig{file=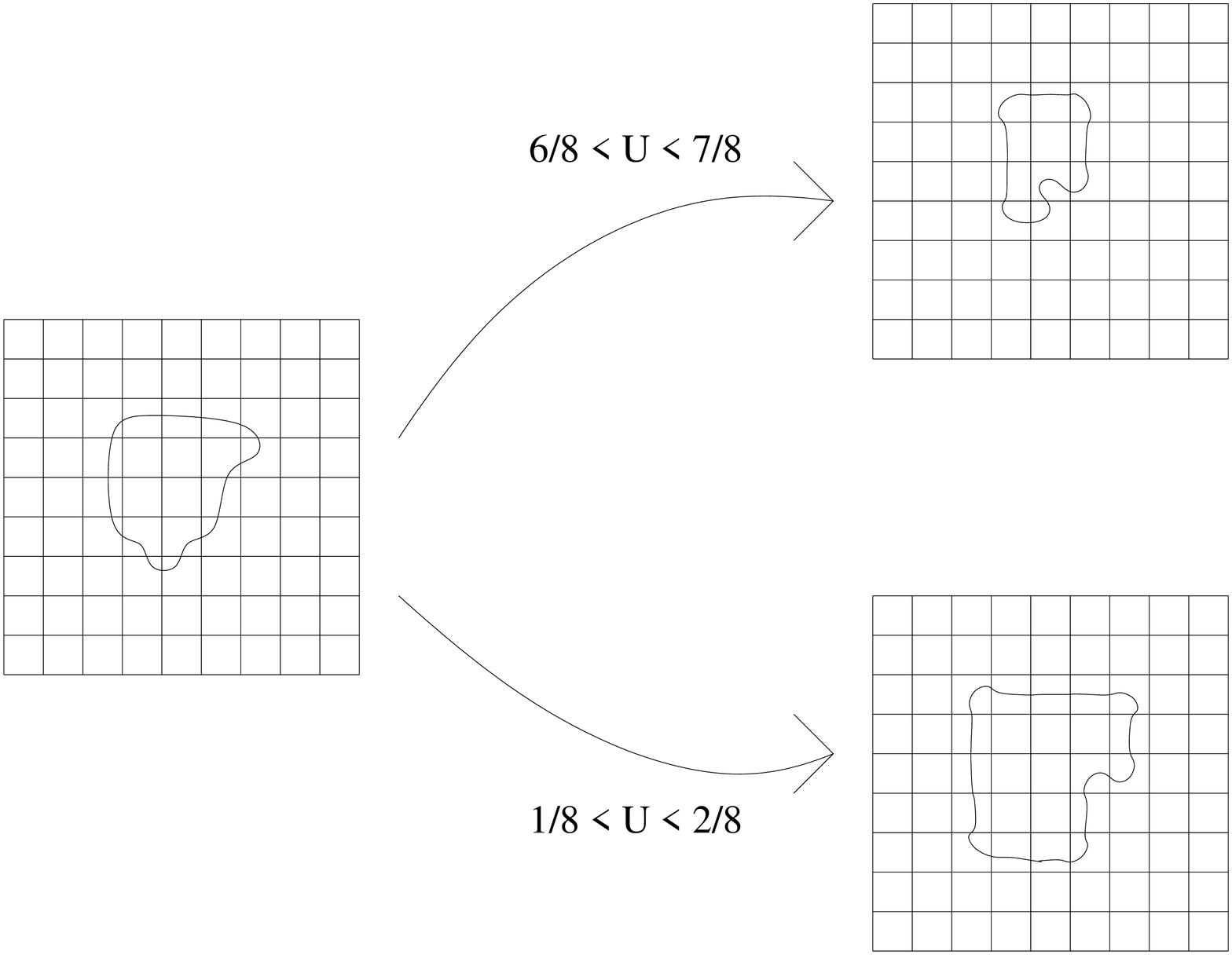, height=3 in, width= 3 in}   
\caption{One step of the evolving set process.}
\end{figure}

\smallskip
\noindent {\bf Definition: Evolving sets.} 
Given $V, \pi$ and $Q$ as above,
consider the Markov chain $\{S_n\}$ on 
{\it subsets} of $V$ with the following transition
rule. If the current state $S_n$ 
is $S \subset V$, choose $U$ uniformly
from $[0,1]$ and let the next state $S_{n+1}$ be
$$
\wS=\{ y: Q(S, y) \geq U \pi(y)\} \,. 
$$
Consequently,
\be \lab{eq:evol}
{\P}(y \in \wS ) = \P\Big(Q(S,y) \geq U \pi(y) \Big)=
\frac{Q(S,y)}{\pi(y)} \,.
\ee
Figure \ref{onestep} illustrates one step of the evolving set process
when the original 
Markov chain is a random walk in a box (with a holding probability of
$\half$).
Since
$\pi$ is the stationary distribution, $\emptyset$ and $V$ are
absorbing states for the evolving set process. 

Write
$\P_S\Big( \cdot \Big):= \P\Big(\,\cdot\, \l| \, S_0 = S \Big)$
and similarly for $\e_S\Big( \cdot \Big)$.
The utility of evolving sets stems from the relation
\[
p^n(x,y) = {\pi(y) \over \pi(x)}\;\Px\left(y \in S_n \right)
\]
(see Proposition \ref{keyprop}). Their connection to mixing is
indicated by the inequality
\[
\Vert \mu_n - \pi \Vert \leq \frac{1}{\pi(x)}
\Ex \sqrt{ \pi(S_n) \wedge \pi(S_n^c}) \, ,
\]
 where $\mu_n:=p^n(x, \, \cdot \,)$; see (\ref{chibound}) for a
 sharper form of this. The connection of evolving sets to
conductance can be seen in Lemma \ref{ll} below.

\begin{figure}
\label{holes}
\centering
\epsfig{file=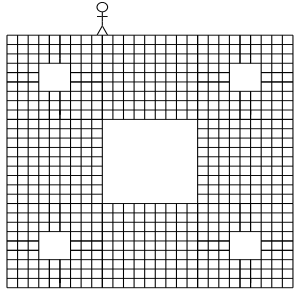, width = 1.5 in, height=1.5 in}
\caption{A box with holes.}
\end{figure}

{\bf Example 1 (Random Walk in a Box): } Consider a simple random walk in
an $n \times n$ box. 
To guarantee condition (\ref{pain}) we add a holding probability of
$\half$ to each state (i.e., with probability $\half$ do nothing, else
move as above). The conductance profile satisfies 
\[
\phi(u) \geq {a \over n \sqrt{u}}
\]
for $1 \leq u \leq 1/2$, where $a$ is a constant. 
Thus our bound implies that the $\epsilon$ uniform mixing time is at most
\[
C_{\epsilon}+4\int_{1/n^2}^{1/2} 
{1 \over u \left( { a \over n\sqrt{u} }
\right)^2}\, du =O(n^2),
\]
which is the correct order of magnitude. Of course, other techniques
such as coupling or spectral methods would give the
correct-order bound of $O(n^2)$ in this case. However, 
these techniques are not robust under small perturbations of the
problem, whereas the conductance method is.\\
\\
{\bf Example 2 (Box with Holes):} 
For a random walk in a box with holes (see Figure \ref{holes}), 
it is considerably harder to apply
coupling or spectral methods. However, it is clear that
the conductance
profile for the random walk is unchanged (up a constant factor), 
and hence the mixing
time  is still $O(n^2)$.

{\bf Example 3 (Random Walk in a Percolation Cluster):}
In fact, the conductance method is robust enough to handle an even
more extreme variant: Suppose that each edge in the box is deleted
with probability $1-p$, where $p > \half$. Then with high
probability there is a connected component that contains a constant
fraction of the original edges. 
Benjamini and Mossel \cite{BM} 
showed that for the random walk in the big
component  the conductance profile is sufficiently close 
(with high probability) 
to that of the box and deduced that the mixing time is still $O(n^2)$.
(See~\cite{MR} for analogous results in higher dimensions.)
By our result, this also applies to the uniform mixing times.
\begin{figure}
\label{perc}
\centering
\epsfig{file=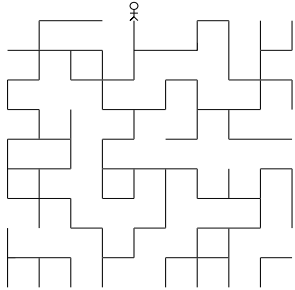, width=1.5 in, height=1.5 in}   
\caption{Random walk in a percolation cluster.}
\smallskip
\end{figure}

{\bf Example 4 (Random Walk on a Lamplighter Group): } 
The following natural chain  mixes more rapidly in the sense of
total variation than in the uniform sense. 
A state of this chain consists of $n$ 
lamps arrayed in a circle,
each lamp either  {\em on\/} (1) or {\em off\/} (0),
 and a
{\em lamplighter\/} located next to one of the lamps.
In one ``active'' step of the chain, the lamplighter either switches the 
current lamp or moves at random to one of the two adjacent lamps.
We consider the lazy chain that stays put with probability $1/2$ and
makes an active step with probability $1/2$.
The path of the lamplighter is a delayed simple random walk on a
cycle, and this implies that $\tau_V(1/4)=\Theta(n^2)$,
see~\cite{HJ}. However, by considering the possibility that
the lamplighter stays in one half of the cycle for a long time,
one easily verifies that $\tau(1/4) \ge c_1 n^3$ for some constant $c_1>0$.
Using the general estimate 
$\tau(\epsilon) =O(\tau_V(\epsilon) \log(1/\pi_*))$
gives a matching upper bound
$\tau(1/4)=O(n^3)$. 

\begin{figure}
\centering
\epsfig{file=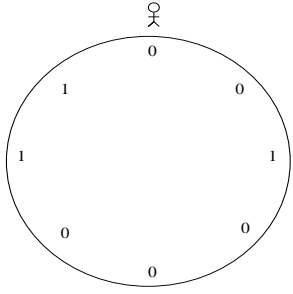, width=1.5 in, height=1.5 in}  
\caption{Random walk on a lamplighter group}
\end{figure}
\smallskip

\section{Further results and proof of Theorem 1} \lab{fur}
We will actually prove a stronger form of Theorem \ref{hk},
using the boundary gauge
$$
\psi(S) := 1 - \e_S \sqrt{ \frac{\pi(\wS)}{\pi(S)}}
$$
instead of the conductance $\phi_S$. 
The next lemma relates these quantities.
\begin{lemma} \label{ll} Let $\emptyset \neq S \subset V$. 
If (\ref{pain}) holds, then $\psi(S) \geq \phi^2_S/2$. 
More generally, if $0<\gamma \le \half$ and
$p(x,x) \ge \gamma$ for all $x \in V$,
 then $\psi(S) \geq \frac{\gamma^2}{2(1-\gamma)^2} \phi^2_S$. 
\end{lemma}
See \S \ref{cp} for the proof. In fact, $\psi(S)$ 
is often much larger than $\phi^2_S$. 

Define the {\em root profile\/} $\psi(r)$ for $r \in [\pim,1/2]$ by
\be \lab{defgam}
\psi(r) = \inf \{ \psi(S): \pi(S) \leq r \},
\ee
and for $r>1/2$, let $\psi(r):=\psi_*=\psi(\half)$.
Observe that the root profile 
$\psi$ is (weakly) decreasing on $[\pim,\infty)$.

For a measure $\mu$ on $V$, write
\be \lab{defchi}
\chi^2(\mu,\pi):=
\sum_{y \in V} 
\pi(y) \Big(\frac{\mu(y)}{\pi(y)}-1\Big)^2 = 
\Big( \sum_{y \in V}  \frac{\mu(y)^2}{\pi(y)}\Big) - 1                \,.
\ee
By Cauchy-Schwarz, 
\begin{equation}
\label{fifteen}
2\Vert \mu-\pi\Vert =
\Big\Vert \frac{\mu(\, \cdot \,)} {\pi(\, \cdot \,)} -1 
\Big\Vert_{L^1(\pi)} \le
\Big\Vert \frac{\mu(\, \cdot \,)} {\pi(\, \cdot \,)} -1\Big\Vert_{L^2(\pi)}
=\chi(\mu,\pi) \, .
\end{equation}

We can now state our key result relating evolving sets to mixing.
\begin{theorem}
\label{psith}
Denote $\mu_n=p^n(x, \, \cdot \,)$. Then $\chi^2(\mu_n, \pi) \leq
\epsilon$
for all
$$ 
n \ge
\int_{4\pi(x)}^{4/\epsilon}
\frac{du}{u \psi(u)} \,.
$$ 
%
\end{theorem}
See \S 5 for the proof.

\medskip
\medskip
\noindent{\bf Derivation of Theorem \ref{hk} from Lemma \ref{ll} and
  Theorem \ref{psith}}:

\smallskip

The {\em  time-reversal\/} of a Markov chain on $V$
 with  stationary distribution
$\pi$ and transition matrix $p(x,y)$,
is another Markov chain with stationary distribution $\pi$,
and transition matrix $\hp(\cdot,\cdot)$ that satisfies
$\pi(y)p(y,z)=\pi(z) \hp(z,y)$ for all $y,z \in V$.
Summing over intermediate states gives
$\pi(z) \hp^m(z,y)=\pi(y)p^m(y,z)$ for all $z,y \in V$ and $m \ge 1$.

\smallskip

Since $p^{n+m}(x,z)=\sum_{y \in V} p^n(x,y)p^m(y,z)$,
stationarity of $\pi$ gives
\be \lab{heat}
p^{n+m}(x,z)-\pi(z)=\sum_{y \in V} \Big(p^n(x,y)-\pi(y)\Big) 
                   \Big(p^m(y,z)-\pi(z) \Big)
\ee
whence 
\begin{eqnarray} \lab{heat2}
  & & \Bigl| \frac{p^{n+m}(x,z)-\pi(z)}{\pi(z)} \Bigr| \\[1ex]
&=&\Bigl| \sum_{y \in V} \pi(y) \Big(\frac{p^n(x,y)}{\pi(y)}-1\Big) 
                   \Big(\frac{\hp^m(z,y)}{\pi(y)} -1 \Big) \Bigr| \\[1ex]
&\le & 
\chi\Big(p^n(x,\cdot), \pi \Big) \chi\Big(\hp^m(z,\cdot), \pi \Big) \label{s}
\end{eqnarray}
by Cauchy-Schwarz.    

The quantity $Q(S,S^c)$ represents,
for any $S \subset V$, the asymptotic 
frequency of transitions from $S$ to $S^c$ in the stationary Markov
chain with transition matrix $p(\cdot,\cdot)$ and hence 
$Q(S,S^c)=Q(S^c,S)$.
It follows that the
time-reversed chain has the same conductance
profile $\phi(\cdot)$ as the original Markov chain. 
Hence, Lemma \ref{ll} and
Theorem \ref{psith} imply that if 
\[
m, \ell \ge
\int_{ 4(\pi(x) \wedge \pi(y)) }^{4/\epsilon}
\frac{2 du}{u \phi^2(u)}, 
\]
and (\ref{pain}) holds, then 
$$
\chi\Big(p^\ell(x,\cdot), \pi \Big) \le \sqrt{\epsilon}
\; \mbox{{ \rm and }} \;  \chi\Big(\hp^m(z,\cdot), \pi \Big) \le
\sqrt{\epsilon} \, .
$$
 Thus by (\ref{s}),
\[
\Bigl|\frac{p^{\ell+m}(x,z)-\pi(z)}{\pi(z)} \Bigr| \leq \epsilon \, ,
\]
and Theorem \ref{hk} is established. 

In fact, the argument above
yields the following more general statement.

\begin{theorem} 
\label{hk2}
Suppose that $0<\gamma \le \half$ and $p(x,x) \ge \gamma$ for all $x \in V$.
If
\be \lab{inteq2}
n \ge 1 +
\frac{(1-\gamma)^2}{\gamma^2} \int_{ 4(\pi(x) \wedge \pi(y)) }^{4/\epsilon}
\frac{4 du}{u \phi^2(u)} \,,
\ee
then  (\ref{unif}) holds.

\end{theorem}

To complete the proof
of Theorems \ref{hk} and \ref{hk2},
it suffices to prove Lemma \ref{ll} and Theorem \ref{psith}.
This is done in \S \ref{cp} and \S \ref{psi}, respectively.

\section{Properties of Evolving Sets} \lab{pro}

\begin{lemma} \lab{lem:mart}
The sequence $\{\pi(S_n)\}_{n \ge 0}$ forms a martingale. 
\end{lemma}
\begin{proof}
By (\ref{eq:evol}), we have 
\begin{eqnarray*}
\e\left( \pi(S_{n+1}) \Bigl| S_n\right) 
&=& \sum_{y \in V} \pi(y)
   \;{\P}\left(y \in S_{n+1} \Bigl| S_n\right) \\[1ex]
                   &=& \sum_{y \in V} \; Q(S_n,y)
                   \, = \, \pi(S_n) \, .
\end{eqnarray*}
\end{proof}
The following proposition relates the
 $n$th order transition probabilities
 of the original chain to the evolving set process. 
\begin{proposition} \label{keyprop}
For all $n \geq 0$ and $x,y \in V$ we have
\[
p^n(x,y) = {\pi(y) \over \pi(x)}\;\Px\left(y \in S_n \right).
\]
\end{proposition}
\begin{proof}
The proof is by induction on $n$. The case $n=0$ is
trivial. Fix $n > 0$ and suppose that the result holds for
$n-1$. Let $U$ be the uniform random variable used to generate $S_n$ from
$S_{n-1}$. Then 
\begin{eqnarray*}
p^n(x,y) &=& \sum_{z \in V} p^{n-1}(x,z)p(z,y) \\[1ex]
&=&  \sum_{z \in V} \Px\left(z \in S_{n-1}\right) 
{\pi(z) \over \pi(x)} \,p(z,y) \\[1ex]
&=&  {\pi(y) \over \pi(x)} \Ex \left({1 \over \pi(y)}Q(S_{n-1}, y)\right) 
 \\[1ex]
& = &   {\pi(y) \over \pi(x)} \;\Px(y \in S_n).
\end{eqnarray*}
\end{proof}

We will also use the following duality property of evolving sets.
\begin{lemma}
\label{comp}
Suppose that $\{S_n\}_{n \ge 0}$ 
is an evolving set process. Then the sequence of complements
$\{S_n^c\}_{n \ge 0}$ 
is also an evolving set process, with the same transition probabilities.
\end{lemma}
\begin{proof}
Fix $n$ and 
let $U$ be the uniform random variable used to generate $S_{n+1}$ from
$S_n$.  Note that $ Q(S_n, y) + Q(S_n^c, y) = Q(V, y) = \pi(y)$. 
Therefore, with probability 1,
\begin{eqnarray*}
S_{n+1}^c &=& \Big\{y: Q( S_n, y) < U \pi(y) \Big\}\\ 
&=& \Big\{y: Q(S_n^c, y) \geq (1-U) \pi(y) \Big\}.
\end{eqnarray*}
Thus, $\{S_n^c\}$ has the same transition probabilities as $\{S_n\}$,
since $1-U$ is uniform. 
\end{proof}

Next, we write the $\chi^2$ distance between $\mu_n:=p^n(x,\, \cdot\,)$
and $\pi$ in terms of
evolving sets. Let $\{S_n\}_{n \ge 0}$ and $\{\Lambda_n\}_{n \ge 0}$
be two independent replicas of the evolving set process, with
$S_0=\Lambda_0=\{x\}$.
Then by (\ref{defchi}) and Proposition \ref{keyprop},
$\chi^2(\mu_n,\pi)$ equals
\begin{eqnarray} 
& &
\sum_{y \in V} 
\pi(y) \frac{\Px(y \in S_n)^2}{\pi(x)^2}-1 \\[1ex]
&=&\frac{1}{\pi(x)^2} \Big[\sum_{y \in V} \pi(y)
\Px\Big(\{y \in S_n\}\cap \{y \in \Lambda_n\}\Big) 
 -\pi(x)^2 \Big] \\[1ex]
&=& 
\label{calcchi}
\frac{1}{\pi(x)^2} 
\Ex\Big(\pi(S_n \cap \Lambda_n)-\pi(S_n)\pi(\Lambda_n) \Big)
\,,
\end{eqnarray}
where the last equation uses the relation
$\pi(x)=\Ex \pi(S_n)=\Ex\pi(\Lambda_n)$.
For any two sets $S,\Lambda \subset V$, 
\[
\pi(S \cap \Lambda) + \pi(S^c \cap \Lambda) = 
\pi(\Lambda) = \pi(S) \pi(\Lambda) + \pi(S^c) \pi(\Lambda),
\]
and hence
\[
|\pi(S \cap \Lambda) - \pi(S)\pi(\Lambda)| =  
|\pi(S^c \cap \Lambda) - \pi(S^c)\pi(\Lambda)|.
\]
Similarly, this expression doesn't change if we replace $\lam $ by
$\lam^c$. 
Thus, if we denote
\[
S^{\sharp}:= 
\left\{\begin{array}{ll}
S & 
\mbox{if $\pi(S) \leq \half$;}\\
S^c & \mbox{otherwise,}\\
\end{array}
\right.
\]
then 
\begin{eqnarray*}
|\pi(S \cap \Lambda) - \pi(S)\pi(\Lambda)| &=&  
|\pi(S^{\sharp} \cap \Lambda^{\sharp}) 
  - \pi(S^\sharp)\pi(\Lambda^\sharp)| \\
 &\le&  
   |\pi(S^\sharp)\wedge \pi(\Lambda^\sharp)| \\
 &\le&
   \sqrt{\pi(S^\sharp)\pi(\Lambda^\sharp)}  \,.
\end{eqnarray*}
Inserting this into (\ref{calcchi}), we obtain
$$
\chi^2(\mu_n,\pi) \le 
\frac{1}{\pi(x)^2}
\e \sqrt{\pi(S^\sharp_n)\pi(\Lambda_n^\sharp)}  \,,
$$
whence 
\be \lab{chibound}
2 \| \mu_n-\pi\| \le \chi(\mu_n,\pi) \le \frac{1}{\pi(x)}
\e \sqrt{\pi(S^\sharp_n)} \, .
\ee 

\section{Evolving sets and conductance profile: proof of Lemma 3}
\label{cp}
\begin{lemma}
\label{comp2}
For every real number $\beta \in [-\half,\half]$, we have
\[
\frac{ \sqrt{1+2\beta} + \sqrt{1-2\beta}}{2} \leq \sqrt{1 - \beta^2}
\leq 1 - \beta^2/2.
\]
\end{lemma}
\begin{proof}
Squaring gives the second inequality and converts the first inequality
into
\[
1+2\beta + 1-2\beta +2 \sqrt{1 - 4 \beta^2} \leq 4 (1 - \beta^2)
\]
or equivalently, after halving both sides,
\[
\sqrt{1 - 4 \beta^2} \leq 1 - 2\beta^2 \,,
\]
which is verified by squaring again.
\end{proof}
\begin{lemma} \lab{lemphi}
Let
\be \lab{defvphi}
\varphi_S:=\frac{1}{2\pi(S)}\sum_{y \in V} \Big(Q(S,y) \wedge Q(S^c,y)
\Big) \,.
\ee 
Then
\be \lab{ine}
1-\psi(S) \leq \frac{ \sqrt{1+2\varphi_S} + \sqrt{1-2\varphi_S}}{2}
\leq 1 - \varphi_S^2/2 \, .
\ee
\end{lemma} 
\begin{proof}
The second inequality in (\ref{ine}) follows immediately 
from Lemma \ref{comp2}. To see the first
inequality, let
$U$ be the uniform random variable used to generate
$\wS$ from $S$. Then
$$
\P_S\Big(y \in \wS \, \l| \, U<\half\Big)= 1 \wedge \frac{2Q(S,y)}{\pi(y)} \,.
$$
Consequently,
$$
\pi(y) \P_S(y \in \wS \, | \, U<\half)= 
Q(S,y)+\Big(Q(S^c,y) \wedge Q(S,y)\Big) \,.
$$
Summing over $y \in V$, we infer that
\be
\e_S \Big( \pi(\wS) \, \l|\, U < \half \Big) = 
\pi(S) + 2 \pi(S)\varphi_S \,.
\ee
Therefore, $\RR:= \pi(\wS)/ \pi(S)$ 
satisfies $\e_S(\RR | U < \half) = 1 + 2 \varphi_S$. 
Since $\e_S \RR = 1$, it follows that 
$$ \e_S(\RR \, | \, U \geq \half) = 1 - 2 \varphi_S \, .
$$
Thus
\begin{eqnarray*}
1-\psi(S) &=& \e(\sqrt{\RR}) \\
&=& 
\frac{ \e(\sqrt{\RR} \bigl|U < \half) +\e(\sqrt{\RR} \bigl| U \geq \half)}{2}\\
&\leq& 
\frac{ \sqrt{  \e(\RR  | U < \half)} 
+ \sqrt{  \e(\RR  | U \ge \half)}}{2},
\end{eqnarray*}
by Jensen's inequality (or by Cauchy-Schwarz). This completes the proof. 
\end{proof}
\begin{proofof}{Proof of Lemma \ref{ll}}
If $p(y,y) \ge 1/2 \; \forall y \in V$,
then it is easy to check directly that 
$\varphi_S =\phi_S$ for all $ S \subset V$.

If we are only given that $p(y,y) \ge \gamma\; \forall y \in V$,
where $0<\gamma  \le \half$, we can still conclude that for $y \in S$,
$$
Q(S,y) \wedge Q(S^c,y) \ge \gamma\pi(y) \wedge Q(S^c,y) \ge
\frac{\gamma}{1-\gamma} Q(S^c,y) \,.
$$
Similarly, for $y \in S^c$ we have 
$Q(S,y) \wedge Q(S^c,y) \ge \frac{\gamma}{1-\gamma} Q(S,y)$.
Therefore 
$$
\sum_{y \in V} [Q(S,y) \wedge Q(S^c,y)] \ge 
\frac{2\gamma}{1-\gamma} Q(S,S^c)\, ,
$$
 whence
$\varphi_S \ge \frac{\gamma}{1-\gamma} \phi_S$.
This inequality, in conjunction with Lemma \ref{lemphi},
yields Lemma \ref{ll}.
\end{proofof}

\section{Proof of Theorem \ref{psith}} \lab{psi}
Denote by $\k(S,A)=\P_S(\wS=A)$ the transition kernel for the evolving set
process. In this section we will use another Markov chain on sets with
transition kernel
\begin{equation}
\label{*}
\kh(S,A)= \frac{\pi(A)}{\pi(S)} \k(S,A).
\end{equation}
This is the Doob transform of $\k(\cdot,\cdot)$.
As pointed out by J. Fill (Lecture at Amer.\ Inst.\ Math.\, 2004),
the process defined by $\kh$ can be identified
with one of the ``strong stationary duals'' constructed in
\cite{DF}.

The martingale property of the evolving set process,
Lemma \ref{lem:mart}, implies that 
$ \sum_{A} \kh(S,A) = 1$ for all $S \subset V$. The chain 
with kernel (\ref{*}) represents
the evolving set process conditioned to absorb in $V$; we will not use
this fact explicitly.

Note that induction from equation (\ref{*}) gives
\[
\kh^n(S,A) = \frac{\pi(A)}{\pi(S)} \k^n(S,A)
\]
for every $n$, 
since 
\begin{eqnarray*}
\kh^{n+1}(S,B) &=& \sum_A \kh^n(S,A) \kh(A,B) \\
&=& \sum_A \frac{\pi(B)}{\pi(S)} \k^n(S,A)\k(A,B) \\ 
&=& 
\frac{\pi(B)}{\pi(S)} \k^{n+1}(S,B)
\end{eqnarray*}
for every $n$ and $B \subset V$. Therefore, for any function $f$,
\begin{equation}
\label{**}
\eh_S f(S_n) = \e_S \left[ \frac{\pi(S_n)}{\pi(S)} f(S_n) \right],
\end{equation}
where we write $\eh$  for the expectation 
when $\{S_n\}$ has transition kernel $\kh$. 
Define
\[
Z_n = \frac{\sqrt{\pi( S_n^{\sharp})}}{ \pi(S_n)} \,,
\]
and note that $\pi(S_n)=Z_n^{-2}$ when 
$Z_n \ge \sqrt{2}$, that is, when $\pi(S_n) \le \half$.
Then by equations (\ref{**}) and (\ref{chibound}), $\chi(\mu_n, \pi)
\leq \eh_{\{x\}} (Z_n)$ and 
\begin{eqnarray} \lab{pref}
\eh \left( \frac{Z_{n+1}}{Z_n} \Big| S_n\right) &=& 
\e \left( \frac{ \pi(S_{n+1})}{\pi(S_n)} \cdot \frac{Z_{n+1}}{Z_n}
 \Big| S_n \right)     \nonumber\\[1ex]
&=& \e \left( 
\frac{\sqrt{\pi( S_{n+1}^{\sharp})}}{\sqrt{\pi( S_{n}^{\sharp})}}
 \Big| S_n \right) \\[1ex]
&\leq& 1-\psi( \pi(S_n)) = 1 - \f_0(Z_n), \label{fzero}
\end{eqnarray}
where 
$\f_0(z) := \psi(1/z^2)$ is nondecreasing. 
(Recall that we defined $\psi(x) = \psistar$ for all real 
numbers $x \geq \half$.)
Let 
$L_0 = Z_0 =
 \pi(x)^{-1/2}$. 
Next, observe that $\eh(\cdot)$ is just the expectation operator
with respect to a modified distribution, so we can apply
Lemma \ref{rr} below, with $\eh$ in place of $\e$.
By part (iii) of that lemma (with $\delta=\sqrt{\epsilon}$),
for all 
\be \lab{integ}
n \geq \int_{\delta}^{L_0} \frac{2 dz}{z \f_0(z/2)} = 
\int_{\delta}^{L_0} \frac{2dz}{z \psi(4/z^2)} \, ,
\ee
we have $\chi(\mu_n, \pi) \leq \eh_{\{x\}} (Z_n) \leq \delta.$ 
The change of variable $u = 4/z^2$ shows the integral (\ref{integ})
equals
\[
\int_{4\pi(x)}^{4/\delta^2} \frac{du}{u \psi(u)} \leq
\int_{4\pi(x)}^{4/\epsilon} \frac{du}{u \psi(u)}.
\]
This establishes Theorem \ref{psith}. 
\begin{lemma} \label{rr}
Let $\f,\f_0:[0, \infty) \to [0,1]$ be increasing functions. 
  Suppose that $\{Z_n\}_{n \ge 0}$ are non-negative random
  variables with $Z_0=L_0$. Denote $L_n=\e(Z_n)$.
\begin{description}
\item{{\bf (i)}} $\,$ If  
$L_n - L_{n+1}  \geq L_n \f(L_n) $ 
for all $n$, then for every $n \ge \int_{\delta}^{L_0} {dz \over z
  \f(z)}$, we have $L_n \leq \delta$.
\item{{\bf (ii)}}  $\,$
 If $\e(Z_{n+1} | Z_n) \leq Z_n (1 - f(Z_n))$ for all $n$
and the function $u \mapsto uf(u)$ is convex on $(0,\infty)$,
then the conclusion of (i) holds.
\item{{\bf (iii)}}  $\,$
 If $\e(Z_{n+1} | Z_n) \leq Z_n (1 - \f_0 (Z_n))$
 for all $n$ and $\f(z)=\f_0(z/2)/2$, then the conclusion of (i) holds.
\end{description}
\end{lemma}
\begin{proof}
\noindent{{\bf (i)}} $\;$
It suffices to show  that for every $n$ we have
\begin{equation}
\label{ind}
\int_{L_n}^{L_0} {dz \over z \, \f(z)} \geq n.
\end{equation}
Note that for all $k \geq 0$ we have
\[
L_{k+1}\leq L_k \Bigl[ 1 - \f(L_k) \Bigr] \leq L_k e^{- \f(L_k)} \, ,
\]
whence
\[
\int_{L_{k+1}}^{L_k} {dz \over z \f(z)} \geq {1 \over
  \f(L_k)} \int_{L_{k+1}}^{L_{k}} {dz \over z} = {1\over
  \f(L_k)} \log {L_k \over L_{k+1}} \geq 1.
\]
Summing this over $k \in \{0,1,\dots, n-1\}$ gives (\ref{ind}).

\noindent{{\bf (ii)}} $\;$ This is immediate from Jensen's inequality
and (i). 

\noindent{{\bf (iii)}}$\;$ Fix $n \ge 0$. We have
\be \lab{cruc}
\e\left( Z_n - Z_{n+1} \right) \geq \e \left[ 2Z_n \f(2Z_n) \right]
\geq L_n \f(L_n) \,,
\ee 
by Lemma \ref{l8} below. This
yields the hypothesis of (i).
\end{proof}
The following simple fact was used in the proof of Lemma \ref{rr}. 
\begin{lemma}
\label{l8}
Suppose that $Z \geq 0$ is a nonnegative random variable and $\f$ is a
nonnegative increasing function.
Then 
\[
\e\Bigl( Z\f(2Z) \Bigr) \geq 
\frac{\e Z}{2} \cdot \f( \e Z).
\]
\end{lemma}
\begin{proof}
Let $A$ be the event $\{Z \geq  \e Z/2 \}$. Then
$\e(Z \one_{A^c}) \leq \e Z/2$,
so $\e(Z \one_A) \geq \e Z/2$. Therefore,
\[
\e\Bigl( Z\f(2Z) \Bigr) \geq 
\e\Bigl( Z \one_A \cdot \f( \e Z) \Bigr) \ge 
\frac{\e Z}{2}  \f(\e Z)  \,.
\] \end{proof}

\section{Infinite stationary measures: proof of Theorem \ref{hki}}
\begin{proof}
For a probability  measure $\mu$ on $V$, define $\chi^2(\mu, \pi)$ by
\be \lab{i:defchi}
\chi^2(\mu,\pi):=
\sum_{y \in V} 
\pi(y) \Big(\frac{\mu(y)}{\pi(y)}\Big)^2 = 
 \sum_{y \in V}  \frac{\mu(y)^2}{\pi(y)}               \,.
\ee

We now write $\chi^2(\mu_n,\pi)$ 
in terms of
evolving sets. Let $\{S_n\}_{n \ge 0}$ and $\{\Lambda_n\}_{n \ge 0}$
be two independent replicas of the evolving set process, with
$S_0=\Lambda_0=\{x\}$.
Then by (\ref{i:defchi}) and Proposition \ref{keyprop},

\begin{eqnarray} \lab{i:calcchi}
\chi^2(\mu_n,\pi)&=&
\sum_{y \in V} 
\pi(y) \frac{\p(y \in S_n)^2}{\pi(x)^2} \\[1ex]
&=&\frac{1}{\pi(x)^2} \Big[\sum_{y \in V} \pi(y)
\p\Big(\{y \in S_n\}\cap \{y \in \Lambda_n\}\Big) \Big] \\[1ex]
&=& \frac{1}{\pi(x)^2} 
\e\Big(\pi(S_n \cap \Lambda_n) \Big) \leq 
\frac{1}{\pi(x)^2} \e\left(\sqrt{\pi(S_n)\pi(\Lambda_n)}\right)  \,.
\end{eqnarray}
whence 
\be \lab{i:chibound}
\chi(\mu_n,\pi) \le \frac{1}{\pi(x)}
\e \sqrt{\pi(S_n)} \, .
\ee 

As in the finite case, if  
$\kh$  is the Doob transform of $\k$ with respect to $\pi$, then
\begin{equation}
\label{i:**}
\eh_S f(S_n) = \e_S \left[ \frac{\pi(S_n)}{\pi(S)} f(S_n) \right].
\end{equation}
Define
\[
Z_n = \frac{1}{\sqrt{ \pi(S_n)}}.
\]
Then by equations (\ref{i:**}) and (\ref{i:chibound}), $\chi(\mu_n,\pi)
\leq \eh_{\{x\}} (Z_n)$ and 
\begin{eqnarray} \lab{i:pref}
\eh \left( \frac{Z_{n+1}}{Z_n} \Big| \; S_n\right) &=& 
\e \left( \frac{ \pi(S_{n+1})}{\pi(S_n)} \cdot \frac{Z_{n+1}}{Z_n}
 \Big| \, S_n \right)     \nonumber\\[1ex]
&=& \e \left( 
\frac{\sqrt{\pi( S_{n+1}^{})}}{\sqrt{\pi( S_{n}^{})}}
 \Bigl| \,S_n \right) \\[1ex]
&\leq& 1-\psi( \pi(S_n)) = 1 - \f_0(Z_n), \nonumber
\end{eqnarray}
where $\f_0(z) = \psi(1/z^2)$
is increasing. Let 
$L_0 = Z_0 =
 \pi(x)^{-1/2}$. By 
Lemma \ref{rr}(iii) above,
for all 
\be \lab{i:integ}
n \geq \int_{\sqrt{\epsilon}}^{L_0} \frac{2 dz}{z \f_0(z/2)} = 
\int_{\sqrt{\epsilon}}^{L_0} \frac{2dz}{z \psi(4/z^2)} \, ,
\ee
we have $\chi(\mu_n, \pi) \leq \eh_{\{x\}} (Z_n) \leq \sqrt{\epsilon}.$ 
The change of variable $u = 4/z^2$ shows the integral (\ref{i:integ})
equals 
\[
\int_{4\pi(x)}^{4/\epsilon} \frac{du}{u \psi(u)}
\leq \frac{ (1-\gamma)^2}{\gamma^2}
\int_{4\pi(x)}^{4/\epsilon} \frac{2du}{ u \phi^2(u)}.
\]
 
Let $\hp$ denote the time-reversal of $p(\cdot, \cdot)$.
Then for all
\[
m, n \ge
\frac{ (1 - \gamma)^2}{\gamma^2} \int_{ 4(\pi(x) \wedge \pi(y)) }^{4/\epsilon}
\frac{2 du}{u \phi^2(u)} 
\]
we have
$$
\chi\Big(p^n(x,\cdot), \pi \Big) \le \sqrt{\epsilon}
\; \mbox{{ \rm and }} \;  \chi\Big(\hp^m(z,\cdot), \pi \Big) \le
\sqrt{\epsilon} \, .
$$
Thus
\begin{eqnarray} \lab{i:heat3}
\Bigl| \frac{p^{n+m}(x,z)}{\pi(z)} \Bigr|
&=&\Bigl| \frac{1}{\pi(z)}\sum_{y \in V} p^n(x,y) p^m(y,z) \Bigr| \\[1ex]
&=&\Bigl| \sum_{y \in V} \pi(y) \Big(\frac{p^n(x,y)}{\pi(y)}\Big)
                   \Big(\frac{\hp^m(z,y)}{\pi(y)}  \Big) \Bigr| \\[1ex]
&\le & 
\chi\Big(p^n(x,\cdot), \pi \Big) \chi\Big(\hp^m(z,\cdot), \pi \Big)
\leq \epsilon,
\end{eqnarray}
where the first inequality is Cauchy-Schwarz.   
This establishes Theorem \ref{hki}. 
\end{proof}

\section{Continuous Time}
In this section we extend our results to continuous-time, finite chains. We
consider the chain $\{ X_t, t \geq 0 \}$ that moves at rate 1
according to $P$, where $P(x,y)$ is a transition kernel on $V$ with 
stationary distribution $\pi$. Let $\phi_P$ be the conductance profile 
for $P$.

\begin{theorem} 
\label{cont1}
Let $X_t$ be a continuous-time, finite chain with 
transtion kernel $h_t=e^{t (P-I)}$. 
Then the $\epsilon$-uniform mixing time satisfies
\be \lab{i:taubound}
\tau(\epsilon) \leq \int_{ 4\pi_* }^{4/\epsilon}
\frac{8 du}{u \phi^2(u)} \,.
\ee
More precisely, 
if
\be \lab{i:inteq}
t \ge
\int_{ 4(\pi(x) \wedge \pi(y)) }^{4/\epsilon}
\frac{8 du}{u \phi^2(u)} \, ,
\ee
then 
\be \lab{i:unif}
\l| \frac{h_t(x,y) - \pi(y)}{\pi(y)} \r| \leq \epsilon.
\ee
\end{theorem}

\begin{proof}
As before, it is enough to show that 
$\chi^2\left(
\P(X_t = \cdot), \pi\right) \leq
\epsilon$
for all
$$ 
t \ge
\int_{4\pi(x)}^{4/\epsilon}
\frac{4du}{u \phi^2(u)} \,.
$$ 

Consider the Markov operator $\lP = \half (P + I)$ with corresponding transition probabilities $\lp(\cdot,\cdot)$.
Let $\lphi$ and $\lpsi$ be the conductance profile and
root profile of $\lP$, respectively. Note  $\lP$
satisfies condition (\ref{pain}) so
Theorems \ref{psith} and \ref{hk} apply.  
Let 
$\{\xt_t\}$ be the chain with transition kernel $\hc_t = e^{t(\lP - I)}$.
Observe that $e^{2t (\lP - I)} = e^{t (P - I)}$, so $\xt_{2t}$
has the same law as $X_{t}$. 
Let 
$\mu_t=\hc_t(x, \, \cdot \,)=h_{t/2}(x, \, \cdot \,)$. 
Since $\lphi = \half \phi$ and $\lpsi \geq {\lphi^2 \over 2}$, 
it is enough to show that
$\chi^2(\mu_t, \pi) \leq
\epsilon$
for all
$$ 
t \ge
\int_{4\pi(x)}^{4/\epsilon}
\frac{du}{u \psi(u)} \,,
$$ 
We accomplish this using the 
natural continuous-time evolving set process
corresponding to $\xt_t$. Let $\{ \Sc(t): t \geq 0 \}$ be 
the process which at rate 1 moves
according to the evolving set transition kernel for $\lP$. 
Let $\{\Sc_n: n \geq 0\}$
be the (discrete time) evolving set process for $\lP$.  
Note that
\begin{eqnarray}
\P_x(\xt_t = y) 
&=& \sum_{j=0}^{\infty} \left[ e^{-t} \; {t^j \over j!} \right] 
\lp^j(x,y) \\ 
&=& \sum_{j=0}^{\infty} \left[ e^{-t} \; {t^j \over j!} \right] 
{\pi(y) \over \pi(x)}\;\Px\left(y \in S_j \right)
 \\ 
&=& 
{\pi(y) \over \pi(x)}\;\Px\left(y \in \Sc(t) \right).
\end{eqnarray}
Our proof will parallel the proof of Theorem \ref{psith}.
One can argue as in Section \ref{pro} to obtain 
$\chi(\mu_t,\pi) \le \frac{1}{\pi(x)}
\e \sqrt{\pi(S^\sharp_t)} \, .$
Define
\[
\zc_t = \frac{\sqrt{\pi( \Sc_t^{\sharp})}}{ \pi(\Sc_t)} \,,
\]
and let $\lc_t = \eh(\zc_t)$, so that $\chi(\mu_t, \pi) \leq \lc_t$. 
Note that 
\[
\lc_t = \sum_{j=0}^{\infty} \left[ e^{-t} \; {t^j \over j!} \right]  L_j
\]
is differentiable. 
Equation (\ref{fzero}) implies that 
\[
\eh \left( \frac{\zc_{t+\epsilon}}{\zc_t} \Big| \; \Sc_t, D \right)  
\leq 1 - \f_0(\zc_t)
\]
where $D$ is the event that the evolving set 
process $\{\Sc(\cdot)\}$ makes exactly one transition
in $(t, t+ \epsilon]$. It follows that for all $t \geq 0$ we have
\[
\e( \zc_{t} - \zc_{t + \epsilon} | \Sc_t) \geq  \epsilon \zc_t f_0(\zc_t)
+ O( \epsilon^2).
\]
Fix $t \ge 0$. Taking expectations above, we get
\be 
\e\left( \zc_t - \zc_{t+\epsilon} \right) \geq \epsilon \e 
\left[ 2\zc_t \f(2\zc_t) \right] + O(\epsilon^2)
\geq \epsilon \lc_t \f(\lc_t) + O(\epsilon^2)\,,
\ee 
where the last inequality holds by Lemma \ref{l8}.

It follows that $\lc_{t+ \epsilon} - \lc_t \leq - \epsilon 
\lc_t f( \lc_t) + O(\epsilon^2)$ and hence
\be
\label{diff}
\lc'_t \leq - \lc_t f( \lc_t).
\ee
The following Lemma is an analog of Lemma \ref{rr}.
\begin{lemma}
\label{clem}
For every 
\[
t \geq \int_{\delta}^{\lc_0} {dz \over z f(z)},
\]
we have $\lc_t \leq \delta$. 
\end{lemma}
\begin{proof}
It's enough to show that for all $t \geq 0$ we have
\be
\label{check}
t \leq \int_{\lc_t}^{\lc_0} { dz \over z f(z) }.
\ee
This is an equality for $t = 0$, and differentiating both sides gives
\[
1 \leq - {1 \over \lc_t f(\lc_t)} \lc'_t,
\]
which holds by equation (\ref{diff}). 
\end{proof}
Lemma \ref{clem} implies that
for all 
\be \lab{i:integc}
t \geq \int_{\delta}^{\lc_0} \frac{2 dz}{z \f_0(z/2)} = 
\int_{\delta}^{\lc_0} \frac{2dz}{z \psi(4/z^2)}
\ee
we have $\chi(\mu_t, \pi) \leq \lc_t 
\leq \delta.$ Let $\epsilon=\delta^2$.
We calculated after (\ref{i:integ}) that the integral on the righ-hand
side of (\ref{i:integc}) equals
\[
\int_{4\pi(x)}^{4/\epsilon} \frac{du}{u \psi(u)}.
\]
This establishes Theorem \ref{cont1}.
\end{proof}

\section{Concluding remarks} \lab{con}
\begin{enumerate}
\item The example of the lamplighter group in the introduction shows
  that $\tau_V(1/4)$, 
the mixing time in total variation on the left-hand side of
(\ref{LK}), can be considerably 
smaller than the corresponding  uniform mixing time $\tau(1/4)$
(so an upper bound for $\tau(\cdot)$ is strictly  stronger).
We note that there are simpler examples of this phenomenon.
For lazy random walk on a clique of $n$ vertices,
$\tau_V(1/4)=\Theta(1)$ while $\tau(1/4)=\Theta(\log n)$.
To see a simple example with bounded degree, consider a graph 
consisting of two expanders of cardinality $n$ and $2^n$,
respectively, joined by a single edge.
In this case 
$\tau_V(1/4)$ is of order $\Theta(n)$, while 
$\tau(1/4)=\Theta(n^2)$. 

\item
\label{spectralgap}
Let $X_n$ be a finite, reversible chain with transition matrix $P$. 
Write $\muxn := p^n(x, \cdot)$. 
Equation (\ref{chibound}) gives
\be 
\lab{chibound2}
\chi(\muxn,\pi) \le \frac{1}{\pi(x)}
\e \sqrt{\pi(S^\sharp_n)}
\le \frac{1}{\sqrt{\pi(x)}}
(1 - \psistar)^n \, .
\ee 

Let $f_2: V \to \R$ be the second eigenfunction of $P$ 
and $\lambda_2$ the second eigenvalue, 
so that 
$Pf_2 = \lambda_2 f_2$.
For $x \in V$, define $f_x: V \to \R$ by $f_x(y) = \delta_x(y) - \pi(y)$,
where $\delta$ is the Dirac delta function. 
We can write
$f_2 = \sum_{x \in V} \alpha_x f_x$. 
Hence
\begin{eqnarray}
\Big\Vert {P^n f(\cdot) \over \pi(\cdot)}\Big\Vert_{L^2(\pi)} 
&\leq& 
\sum_x \alpha_x 
\Big\Vert {P^n f_x(\cdot) \over \pi(\cdot)}\Big\Vert_{L^2(\pi)} \\
&=&
\sum_x \alpha_x \chi(\muxn, \pi) \\
&\leq& {\rm const}\cdot \max_x \chi(\muxn, \pi) \\
\label{eI}
&\leq& {\rm const}\cdot (1 - \psistar)^n,
\end{eqnarray}
where the first line is subadditivity of a norm and the last line
follows from (\ref{chibound2}). 
But
\be
\lab{eII}
\Big\Vert {P^n f(\cdot) \over \pi(\cdot)}\Big\Vert_{L^2(\pi)} 
\geq
\Big\Vert {P^n f(\cdot) \over \pi(\cdot)}\Big\Vert_{L^1(\pi)} 
=
\sum_x |P^n f_2(x)| = \lambda_2^n \sum_x |f_2(x)|.
\ee
Combining (\ref{eI}) and (\ref{eII}) gives
$\lambda_2^n \leq c\cdot (1 - \psistar)^n$ for a constant $c$.
Since this is true for all $n$, we must have $\lambda_2 \leq 
1 - \psistar$, so $\psistar$ is a lower bound for the spectral gap.

\item 

Variants of conductance can give better bounds on $\psi_S$. 
For $S \subset V$, define
\[
\theta_S:=\frac{1}{\pi(S)} \sum_{y \in S} \sqrt{ \pi(y) Q(S^c,y)}.
\]
Note that for reversible chains 
we have $$\theta_S = 
\frac{1}{\pi(S)}\sum_{y \in S} \pi(y) \sqrt{p(y, S^c)},$$ which is 
strictly greater than $\phi_S$    for
$S \notin \{ \emptyset, V\}$ (since $\phi_S$ can be written in a
similar way, but without the square root.)

Following Houdr\'e and Tetali~\cite{HT}, denote
\be \lab{h2+}
h_2^+ = \inf \left\{  \theta_S : \pi(S) \leq \half \right\} \, .
\ee
\begin{theorem} Suppose that $p(x,x) \geq
  \half$ for all $x$. Then
\[
\psi(S) \geq \frac{ \theta_S^2} { 8 \log(2/\theta_S^2)} \, .
\]
Consequently, assuming reversibility, the spectral gap $1-\lambda_2$
satisfies 
\be \lab{eqht}
1-\lambda_2 \ge  \frac{ (h_2^+)^2} { 8 \log\Big(2/(h_2^+)^2\Big)}
 \ge  c \, \frac{ (h_2^+)^2} { |\log (h_2^+)| } \,.
\ee
\end{theorem}
Up to the logarithmic factor in the denominator,
this type of inequality was conjectured by  
Houdr\'e and Tetali~\cite[Remark 3.5]{HT}.
\begin{proof}
For $u \in [0,1]$, let $A_u = \{y: Q(S,y) > u \pi(y) \}.$ Note that
$A_u \subset S$ for $u > \half$ since 
$p(x,x) \geq \half$ for all $x$. We have
\begin{eqnarray*}
\e \sqrt{ { \pi(\wS) \over \pi(S)}} &=& 
\int_0^1 \sqrt{ { \pi(A_u) \over \pi(S)}} du \\
&=& \int_0^1 \sqrt{1 +  { \pi(A_u) - \pi(S) \over \pi(S)}} du \\
&\leq& 1 + \half \int_0^1 { \pi(A_u) - \pi(S) 
\over \pi(S)} du \\
& & \hspace{.1 in} 
- \eighth \int_{\half}^1 {(\pi(A_u) 
- \pi(S))^2 \over \pi(S)^2} du, \label{te}
\end{eqnarray*}
by the inequality $\sqrt{1 + t} \leq 1 + t/2 - \eighth t^2 \one(t \leq
0)$, valid for $t \geq -1$. 
Define
\begin{eqnarray*}
B_t := S - A_{1-t} &=& \{y \in S: Q(S,y) \leq (1-t) \pi(y) \}\\
&=& \{y \in S: 
Q(S^c,y) \geq t \pi(y) \}.
\end{eqnarray*}
The middle term in (\ref{te}) vanishes 
by the martingale
property. 
Thus
\begin{eqnarray*}
\psi(S) &\geq& \eighth \int_\half^1 { \pi(S - A_u)^2 \over \pi(S)^2} du \\
&=& \eighth \int_0^\half {\pi(B_t)^2 \over \pi(S)^2} \,dt,
\end{eqnarray*}
where we have made the substitution $t = 1-u$.
Therefore, for any $\alpha \in (0, \half]$ we have
\begin{eqnarray}
8 \pi(S)^2 \psi(S) \cdot \log(1/2\alpha) 
&\geq& \int_{\alpha}^\half \pi(B_t)^2 \,dt \cdot \int_\alpha^\half {1 \over t}
\,dt\\
&\geq& 
\left[ \int_\alpha^\half { \pi(B_t) \over \sqrt{t} } \,dt
\right]^2, \label{*1}
\end{eqnarray}
by Cauchy-Schwarz.
But
\begin{eqnarray*}
\int_\alpha^\half { \pi(B_t) \over \sqrt{t} } \,dt
&=& \sum_{y \in S} \pi(y) \int_\alpha^\half \one(y \in B_t)
{1 \over \sqrt{t}} \,dt \\
&\geq& \sum_{y \in S} \pi(y) \int_\alpha^{Q(S^c,y)/\pi(y)} 
{1 \over \sqrt{t}} \,dt \\
&=& 2 \left( \sum_{y \in S} \sqrt{ Q(S^c,y) \pi(y)} \right)
 - 2 \sqrt{\alpha} \cdot \pi(S) \\
&=& 2 \pi(S) \left( \theta_S - \sqrt{\alpha} \right). 
\end{eqnarray*}
Setting $ \alpha = \theta_S^2/4$ and using equation (\ref{*1}), we get 
\[
8\psi(S) \log(2/\theta_S^2) \geq {1 \over \pi(S)^2}
\left[ \int_{\theta_S^2/4}^\half { \pi(B_t) \over \sqrt{t}} \,dt
\right]^2
\geq \theta_S^2,
\]
and the theorem follows. 
\end{proof}

\item
Theorems \ref{psith} and \ref{hk2} can be improved
under  a convexity condition that holds in many examples.
In the setting of Theorem \ref{psith},
if $\psi(r) \ge \psi_c(r)$ for all $r \ge \pim$
where $u \mapsto u \psi_c(u^{-2})$ is a convex function
of $u \in [0, \infty)$, then $\chi^2(\mu_n, \pi) \leq
\epsilon$ for all  
$$
n \ge \half \int_{ \pi(x) }^{\frac{1}{\epsilon}}
\frac{du}{u \psi_c(u)} \,.
$$
To prove this, follow the proof of Theorem \ref{psith}
until (\ref{pref}), which implies that
$$ 
\eh_n \left( \frac{Z_{n+1}}{Z_n} \Big| S_n \right) 
\le 1-\psi_c( \pi(S_n))
 =1-f(Z_n) \,,
$$
where $f(z):=\psi_c(z^{-2})$. Then invoke Lemma \ref{rr}(ii) 
and apply the change of variable $u=1/z^2$ to the integral there.

Similarly, the following variant of Theorem \ref{hk2} holds.
Suppose that  $p(x,x) \ge \gamma$ for all $x \in V$.
If $\phi(r) \ge \phi_c(r)$ for all $r>0$,
where $u \mapsto u \phi_c^2(u^{-2})$ is a convex function
of $u \in [0, \infty)$, then (\ref{unif}) holds for all
\be \lab{inteq3}
n \ge
\frac{(1-\gamma)^2}{\gamma^2} \int_{ \pi(x) \wedge \pi(y) }^{1/\epsilon}
\frac{2 du}{u \phi_c^2(u)} \,.
\ee
\item
Let ${\mathcal E}$ denote the support
of the evolving set process. 
Theorem \ref{psith} can be improved by using
$\psi_{{\mathcal E}} (r) = 
\inf \{ \psi(S): \pi(S) \leq r, S \in {\mathcal E}\}$,
instead of $\psi(r)$.
For random walk on the $n$ dimensional hypercube,
${\mathcal E}$ consists of Hamming balls,
$\psi_{{\mathcal E}} (r) \ge \frac{c}{n}\log(1/r)$
and
this gives an upper bound of $O( n \log n)$ for 
the uniform mixing time $\tau(1/4)$.

More generally, for any Markov chain $\{X_n\}$
on a poset with a monotone time-reversal,
if $X_0$ is a maximal (or minimal) state, 
then ${\mathcal E}$ consists of increasing (respectively, decreasing)
sets.
\end{enumerate}
\noindent {\bf Acknowledgments.} We are grateful to D. Aldous, L. \lov,
R. Lyons, 
R. Montenegro, E. Mossel and A. Sinclair for useful discussions and comments.

\end{document}